\documentclass[11pt,dvipsnames,table]{article}  

\usepackage[T1]{fontenc}
\usepackage[sc]{mathpazo}

\usepackage[includefoot]{geometry}
\usepackage{amsmath}
\usepackage{amssymb}
\usepackage{amsthm}
\usepackage{array}

\usepackage{tikz}
\usepackage{tikz-3dplot}
\usepackage{pgfplots}
\pgfplotsset{compat=1.15}
\usepackage{tkz-graph}
\usetikzlibrary{arrows}
\usetikzlibrary{arrows.meta}
\usetikzlibrary{automata}
\usetikzlibrary{decorations.markings}
\usetikzlibrary{patterns}
\usetikzlibrary{positioning}
\usetikzlibrary{shapes.geometric}

\usepackage{pdfpages}
\usepackage{microtype}
\usepackage{pifont}
\usepackage{forest}

\usepackage{multirow,multicol,enumitem}
\usepackage{icomma} 
\usepackage{ytableau}
\usepackage{bm}
\usepackage{bbm}
\usepackage{hyperref}
\hypersetup{
    colorlinks=true,
    linkcolor=magenta,
    filecolor=blue,      
    urlcolor=blue,
    citecolor=teal
    }
\usepackage{caption}
\usepackage{mathrsfs}

\theoremstyle{plain}
\newtheorem{theorem}{Theorem}

\theoremstyle{definition}

\theoremstyle{remark}

\usepackage{titlesec}
\titleformat{\section}{\Large\sc\linespread{0.75}\selectfont\hyphenchar\font=-1\relax}{}{0em}%
{}[\vspace{-8pt}\rule{1.0\textwidth}{1pt}\vspace{-10pt}]
\titleformat{\subsection}{\large\bf}{}{0em}{}
\titlespacing{\section}{0em}{2em}{0em}

\newcommand{\name}[1]{\noindent {\it #1}\medskip}
\newcommand{\adr}[1]{\hfill {#1}}
\newcommand{\jointwork}[1]{\noindent {\it This talk is based on joint work with #1}\medskip}
\usepackage[font=small,skip=5pt]{caption}
\usepackage{subcaption}
\usepackage{verbatim}
\usepackage{makecell}
\usepackage{stackrel}
\usepackage{algorithm}
\usepackage[noend]{algpseudocode}
\usepackage[framemethod=tikz]{mdframed}
\algrenewcommand\algorithmicprocedure{\textbf{def}}
\algnewcommand{\LineComment}[1]{\State \(\triangleright\) #1}
\algrenewcommand\algorithmicindent{1.0em} 

\newlength{\bibitemsep}
\newlength{\bibparskip}
\let\oldthebibliography\thebibliography
\renewcommand\thebibliography[1]{%
  \oldthebibliography{#1}%
  \setlength{\parskip}{\bibitemsep}%
  \setlength{\itemsep}{\bibparskip}%
}

\surroundwithmdframed[linecolor=black,linewidth=1pt,skipabove=2mm,skipbelow=2mm,leftmargin=0mm,rightmargin=0mm,innerleftmargin=2mm,innerrightmargin=2mm,innertopmargin=1.7mm,innerbottommargin=1mm]{algo}

\newenvironment{algo}[2]%
{{\bfseries #1} {\itshape (#2)}.}{}

\newcommand{\avoidSymbol}{\mathsf{Av}}
\newcommand{\avoid}[3][]{\avoidSymbol^{#1}_{#2}(#3)}

\newcommand{\parent}[2][]{p^{#1}(#2)}
\newcommand{\children}[2][]{c^{#1}(#2)}
\newcommand{\childrenRight}[2][]{\overrightarrow{c^{#1}}(#2)}
\newcommand{\childrenLeft}[2][]{\overleftarrow{c^{#1}}(#2)}
\newcommand{\childrenBoth}[2][]{\overleftrightarrow{c^{#1}}(#2)}

\newcommand{\one}{\color{darkblue}1\color{black}}
\newcommand{\two}{\color{darkgreen}2\color{black}}
\newcommand{\thr}{\color{darkred}3\color{black}}
\newcommand{\sml}[1]{\color{black}\underline{#1}\color{black}}
\newcommand{\lrg}[1]{\color{black}\overline{#1}\color{black}}

\definecolor{darkred}{rgb}{0.55, 0.0, 0.0}
\definecolor{darkgreen}{rgb}{0.0, 0.2, 0.13}
\definecolor{darkblue}{rgb}{0.0, 0.0, 0.55}

\begin{document}

    \section*{Exhaustive Generation of Pattern-Avoiding $s$-Words}
  \name{Aaron Williams}
    \adr{Williams College}

\jointwork{Samuel Buick,
Madeleine Goertz,
Amos Lastmann,
Kunal Pal,
Helen Qian,
Sam Tacheny,
Leah Williams, and
Yulin Zhai (NSF Grant~DMS2241623)
}
\vspace{-0.3em}


We introduce a simple approach for generating Gray codes of pattern-avoiding $s$-words (i.e., multiset permutations) and corresponding combinatorial objects.
It generalizes plain changes and the recent \emph{Combinatorial Generation via Permutation Language} series.

\textbf{Introduction}
\emph{Plain changes} is a swap Gray code for $S_n$ the permutations of $\mbox{[n] = \{1,\ldots,n\}}$:
e.g., $
\one\sml{\two}\lrg{\thr}, 
\sml{\one}\lrg{\thr}\two, 
\thr\sml{\one}\lrg{\two}, 
\lrg{\thr}\sml{\two}\one, 
\two\lrg{\thr}\sml{\one}, 
\two\one\thr 
$ 
where a \emph{swap} moves a $\overline{\text{large}}$ digit past a \underline{small} digit.
It was discovered in the 1600s and is also known as the \emph{Steinhaus-Johnson-Trotter algorithm}.
More recently, it has been viewed as a greedy algorithm: ``swap the largest value'' \cite{williams2013greedy}.

A \emph{permutation language} is any $L \subseteq S_n$.
A \emph{jump} moves a $\overline{\text{larger}}$ digit past \textbf{one or~more} \underline{smaller} digits. 
It is \emph{minimal} for $w \in L$ if its \emph{distance} $d$ (i.e., number of smaller digits) is minimized to create $w' \in L$.
Algorithm J debuted at \emph{Permutation Patterns 2019}: ``minimal jump the largest value'' \cite{hartung2022combinatorial}.
It generates zig-zag~languages e.g.,
$
\one\sml{\two}\lrg{\thr}, 
\sml{\one}\lrg{\thr}\two, 
\thr\sml{\one}\lrg{\two}, 
\lrg{\thr}\sml{\two\one}, 
\two\one\thr 
$
for $\avoid{3}{231}$.
But jumps become limited for \emph{$s$-words} (or \emph{$s$-permutations}) which have $s_i$ copies of $i$ for $i \in [m]$.
For example, 
$w = \color{darkblue}1\color{darkgreen}2\color{darkred}333\color{darkgreen}2\color{black}$ is a \emph{Stirling $s$-word} for $s=(1,2,3)$ (i.e., $L = \avoid{s}{212}$) but any jump in $w$ gives an \emph{invalid} $w' \notin L$.
The set of all $s$-words is $S_s$.
Each $v \in [m]$ is a \emph{value}, each copy of a value in $w$ is a \emph{digit}, and a digit's \emph{rank} is by increasing value then left-to-right index.
For example, the ranks of $\color{darkblue}1\color{darkgreen}2\color{darkred}333\color{darkgreen}2\color{black}$ are $\color{darkblue}1\color{darkgreen}2\color{darkred}456\color{darkgreen}3\color{black}$.

We consider ``bumps'' which move a \textbf{run} of some larger value.
Algorithm~B generates many \emph{$s$-word languages} (i.e., subsets of $S_s$).
Several applications are listed below.
\begin{enumerate}[noitemsep,topsep=-\parskip] 
    \small
    \item[(a)] Gray codes for $S_s$ (e.g., $\one\sml{\one}\lrg{\two\two},\one\two\lrg{\two}\sml{\one},\sml{\one}\lrg{\two}\one\two,\two\one\sml{\one}\lrg{\two},\two\sml{\one}\lrg{\two}\one,\two\two\one\one$ for $s=(2,2)$) using transpositions.
    \item[(b)] \emph{Stirling changes} generalizes plain changes to \emph{Stirling $s$-words} $\avoid{s}{212}$ using transpositions.
    \item[(c)] Bump Gray code for regular words counted by $k$-Catalan numbers $\avoid[k-1]{m}{132,121}$ \cite{kuba2012enumeration,williams2023pattern}.
\end{enumerate}
\medskip

Our Gray codes lead to efficient algorithms and have various applications (see Figures).
For example, we generate (b) \emph{looplessly} (i.e., worst-case $O(1)$-time per~word) in Algorithm \ref{alg:fastStirling}.
In addition, (b) leads to a $\pm 1$ tree-inversion Gray code for $s$-increasing trees that proves Theorem \ref{thm:spermutahedronPath}.
Finally, (c) gives Gray codes for various $k$-Catalan objects.

\begin{theorem} \label{thm:spermutahedronPath}
    Every $s$-permutahedron has a Hamilton path. 
    (See \cite{ceballos2019s} and its decreasing trees.)
\end{theorem}
\vspace{-1.0em}  

\begin{figure}[h]    \centerline{\includegraphics[page=4,width=1.0\linewidth]{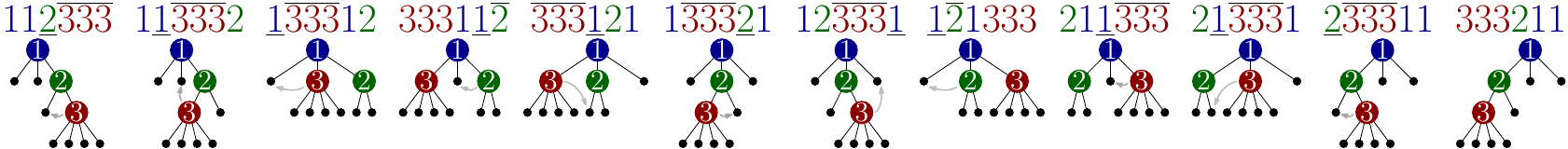}}
    \vspace{-0.055em}
    \caption{
    Stirling changes for Stirling $s$-words with $s=(2,1,3)$ (i.e., permutations of $\{1,1,2,3,3,3\}$ avoiding~$212$) as generated by Algorithm B with its corresponding $s$-increasing tree Gray code.
    }
    \label{fig:sStirling}
    \smallskip
    \centerline{\includegraphics[page=3,width=1.0\linewidth]{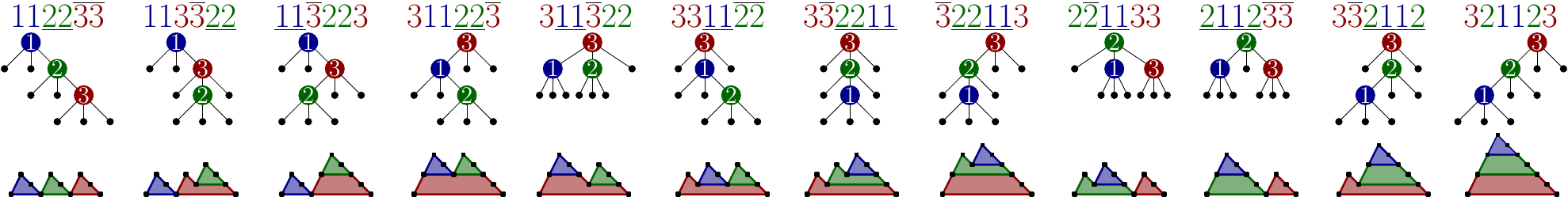}} 
    \vspace{-0.055em}
    \caption{
    Pattern-avoiding $s$-words $\avoid{s}{132,121}$ for $s=(2,2,2)$ with Gray codes for $3$-Catalan objects. 
    The ternary trees differ by moves preserving inorder (i.e., visit self before last~child).
    }
    \label{fig:sCatalan}
    \vspace{-1em}
\end{figure}




\begin{figure}[t]
\centerline{\includegraphics[width=1.1\textwidth]{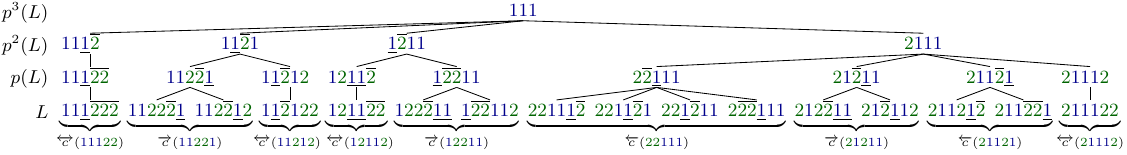}}
    \vspace{-0.25em}
    \caption{Algorithm $B$'s Gray code for $L=\avoid{3,3}{12121}$ and its ancestor languages. 
    Successive children jump the rightmost largest value
    while a bump changes the last child of one parent to the first child of the next parent.
    Note that the rightmost largest value may join these bumps.
    For example, the jump $\two\lrg{\two}\sml{\one}\one\one = \two\one\two\one\one$ in $\parent{L}$ widens to the bump $\two\lrg{\two\two}\sml{\one}\one\one = \two\one\two\two\one\one$ in $L$.
    }
    \label{fig:avoid12121}
    \vspace{-1.1em}
\end{figure}

\textbf{Algorithm B for Bumps}
Let $w = w_1 w_2 \cdots w_n$ be an $s$-word with $s=(s_1, s_2, \ldots, s_m)$.
The \emph{right-maximal run} (or \emph{right-run}) at index $i$ is $w_i w_{i+1} \cdots w_j$ with $w_i = \cdots = w_j$ and $j = n$ or $w_{j} \neq w_{j+1}$. 
A right-run is also determined by $w_i$'s rank.
A \emph{right-bump} moves a right-run to the right past some smaller digits.
Left notions are similar.
A bump's \emph{rank $r$}, \emph{width} $w$, and \emph{distance} $d$ are the rank of $w_i$, $\#$ of larger digits, and $\#$ of smaller digits.
For example, 
$\color{darkblue}1\color{black}\underline{\color{darkblue}1\color{darkgreen}2\color{darkblue}11\color{black}}\color{black}\overline{\color{darkred}33}\color{darkred}3\color{darkblue}11\color{black}$ 
creates~$\color{darkblue}1\color{black}\overline{\color{darkred}33}\underline{\color{darkblue}1\color{darkgreen}2\color{darkblue}11\color{black}}\color{black}\color{darkred}3\color{darkblue}11\color{black}$ 
by a left-bump of $r=9$, $w=2$, and $d=4$.
Bumps can be jumps ($w=1$), transpositions ($d=1$), or swaps ($w=d=1$).

Let $L$ be an $s$-word language.
Suppose a bump at rank $r$ changes $w \in L$ to $w'$.
The bump is \emph{minimal} for its rank and direction if its distance is positive and minimized with $w' \in L$.
Minimal bumps on $w \in L$ are uniquely determined by rank and direction.
Our generalization of Algorithm J \cite{hartung2022combinatorial} replaces values with ranks, and prefers rightwards.

\begin{algo}{Algorithm~$B$}{Greedy Bumps}
This algorithm attempts to generate an $s$-word language $L$ with $s=(s_1, s_2, \ldots, s_m)$ starting from a given or default initial word $w \in L$.
\begin{enumerate}[label={\bfseries B\arabic*.}, leftmargin=7mm, noitemsep, topsep=-1.0em]
\item{} [Initialize] Visit the given initial word $w$ (or by default visit $w = 1^{s_1} 2^{s_2} \cdots m^{s_m}$).
\item{} [Greedy] 
Let $w$ be the most recent word visited.
Visit a new word by applying to $w$ a minimal bump prioritized by largest rank then rightward over leftward. 
Halt if no such bump exists.
Otherwise,~repeat~B2.
\end{enumerate}
\end{algo}
\vspace{-1.6em}




A jump is \emph{maximum} if it uses the longest distance from an index in a direction.
A \emph{zig-zag language} is an $s$-word language closed under maximum jumps (c.f., \cite{hartung2022combinatorial}).
This includes $\avoid{s}{\alpha}$ when $\alpha$'s largest values are \emph{internal} (i.e., not first or last) and \emph{isolated} (i.e., not consecutive).
Zig-zag languages are closed under intersection (and union) and note that the \emph{peakless $s$-words} $\avoid{s}{132,231,121}$ are a subset of every zig-zag language.
%

\begin{theorem} \label{thm:algB}
\!\!\!\footnote{The full paper generalizes this to \emph{zig/zag languages} (rather than zig-zag languages) including $\avoid{s}{212}$.} 
    If $L$ is zig-zag language of $s$-words for $s=(s_1,s_2,\ldots,s_m$), then Algorithm~$B$ generates a bump Gray code for $L$ starting from the non-decreasing word $w = 1^{s_1} 2^{s_2} \cdots m^{s_m} \in~L$.
\end{theorem}
\vspace{-0.85em}  
\vspace{-1.1em} 
\begin{proof}[Sketch]
Consider the end of an inductive argument on $n=\sum s$. 
The \emph{parent} $\parent{s}$ of $s$ is $(s_1,s_2,\ldots,s_m{-}1)$ if $s_m > 1$ or $(s_1,s_2,\ldots,s_{m-1})$ if $s_m = 1$.
Similarly, $\parent{w}$ removes the rightmost $m$ from $w \in L$, and $\parent{L} = \{\parent{w} \mid w \in L\}$. 
The \emph{children} of $w' \in \parent{L}$ are $\children{w'} = \{w \in L \mid \parent{w} = w'\}$.
Since $L$ is a zig-zag language, $\children{w'}$ has $s$-words where the rightmost $m$ is at (a) index $n$, and (b) index $1$ if $s_m=1$ or beside the rightmost $m$ in $w'$ if $s_m>1$; these extremes are equal when $w'$ ends in $m$. 
Let $\childrenRight{w'}$ list $\children{w'}$ in lexicographic order (i.e., the rightmost $m$ jumps left-to-right) and $\childrenLeft{w'}$ in reverse; use $\childrenBoth{w'}$ when $w'$ has one child.
Note that $\parent{L}$ is a zig-zag language, so Algorithm $B$ generates $w'$.
We claim that Algorithm $B$ generates $\childrenRight{w'}$ or $\childrenLeft{w'}$ (or $\childrenBoth{w'}$) when run on $L$.
In other words, it unfurls a \emph{local recursion}.
In particular, a bump transforms the last child of each parent to the first child of the next parent.
See Figure~\ref{fig:avoid12121}.
\end{proof}




\newpage
\clearpage 

\begin{algorithm}[h]
\caption{Loopless generation of Stirling $s$-words $\avoid{s}{212}$ in Stirling changes order.
}
\label{alg:fastStirling}
\hspace{-1em}
\begin{minipage}[t]{0.53\textwidth}
\small
\begin{algorithmic}[1]
\Statex \hspace{-1.5em} $\mathtt{\mathbf{def}\ Stirling(s_1, s_2, \ldots, s_m)}$ \Comment{$1$-based indexing}   \vspace{-0.055em} 
\State $\mathtt{t_1, t_2, t_3, {.}{.}{.}, t_m \leftarrow 0, s_1, s_1{+}s_2, {.}{.}{.}, s_1{+}{\cdot}{\cdot}{\cdot}{+}s_{m-1}}$  
\State $\mathtt{perm \leftarrow 1^{s_1} 2^{s_2} \cdots m^{s_m}}$ \Comment{$s$-word ($s$-perm)}   \vspace{-0.055em}
\State $\mathtt{left \leftarrow [t_1{+}1, {.}{.}{.}, t_{m-1}{+}1]}$ \Comment{left indices~$\mathtt{1,{.}{.}{.},m}$}
\State $\mathtt{inv \leftarrow 0^{m}}$ \Comment{$\#$ of smaller digits right of $\mathtt{1, 2, {.}{.}{.},m}$}
\State $\mathtt{fs \leftarrow 1 2 {\cdot}{\cdot}{\cdot} m}$ \Comment{focus pointers for $\mathtt{1,{.}{.}{.},m}$ \cite{pilau2025efficient,qiu2024generating}}
\State $\mathtt{dirs \leftarrow -1^{m}}$ \Comment{bump directions for $\mathtt{1,{.}{.}{.},m}$}  \vspace{-0.055em}
\State $\mathtt{v \leftarrow \mathtt{fs[m]}}$ \Comment{larger value in next bump}  \vspace{-0.055em}
\While{$\mathtt{v > 1}$} \Comment{$\mathtt{1}$ is never the larger value}  \vspace{-0.055em}
    \State $\mathtt{d} \leftarrow \mathtt{dirs[v]}$ \Comment{direction of $\mathtt{v}$}  \vspace{-0.055em}
    \If{$\mathtt{d} = 1$} \Comment{right-bump $\mathtt{v}$'s run?}  \vspace{-0.055em}
        \State $\mathtt{i \leftarrow left[v]}$  \Comment{first index of $\mathtt{v}$}  \vspace{-0.055em}
        \State $\mathtt{j \leftarrow left[v]+s[v]}$ \Comment{index after $\mathtt{v}$'s run}  \vspace{-0.055em}
    \Else \Comment{left-bump $\mathtt{v}$'s run}  \vspace{-0.055em}
        \State $\mathtt{i \leftarrow left[v]+s[v]-1}$ \Comment{last index of $\mathtt{v}$}  \vspace{-0.055em}
        \State $\mathtt{j \leftarrow left[v]-1}$ \Comment{index before $\mathtt{v}$'s run}  \vspace{-0.055em}
    \EndIf  \vspace{-0.055em}
    \vspace{-0.2em} 
    \State $\mathtt{u \leftarrow perm[j]}$ \Comment{smaller digit to bump}  \vspace{-0.055em}
    \State \textbf{visit} $\mathtt{perm}$  \Comment{next $s$-word in Gray code}  \vspace{-0.055em}
    \State $\mathtt{perm[i] \leftarrow u}$ \Comment{apply bump to the $s$-word}  \vspace{-0.055em}
    \State $\mathtt{perm[j] \leftarrow v}$ \Comment{(as a transposition)}  \vspace{-0.055em}
    \State $\mathtt{left[v] \leftarrow left[v]{+}d}$ \Comment{leftmost $\mathtt{v}$ moved}  \vspace{-0.055em}
    \If{$\mathtt{left[u] = j}$} \Comment{leftmost $\mathtt{u}$ moved?}  \vspace{-0.055em}
        \State $\mathtt{left[u] \leftarrow left[u]{-}d {\cdot} s_v}$ \Comment{$\mathtt{u}$ passed all $\mathtt{v}$}  \vspace{-0.055em}
    \EndIf  \vspace{-0.055em}
    \vspace{-0.2em} 
    \State $\mathtt{inv[v] \leftarrow inv[v]{-}d}$ \Comment{$\mathtt{v}$'s run passed one $\mathtt{u}$}  \vspace{-0.055em}
    \If{$\mathtt{inv[v] = 0}$ \textbf{or} $\mathtt{inv[v] = t_{v}}$} \Comment{$\mathtt{v}$ limit?}  \vspace{-0.055em}
        \State $\mathtt{dirs[v]} \leftarrow -\mathtt{d}$ \Comment{change $\mathtt{v}$'s direction}  \vspace{-0.055em}
        \State $\mathtt{fs[v]} \leftarrow \mathtt{fs[v{-}1]}$  \Comment{inherit focus of $\mathtt{v{-}1}$}  \vspace{-0.055em}
        \State $\mathtt{fs[v{-}1]} \leftarrow \mathtt{v{-}1}$ \Comment{reset focus of $\mathtt{v{-}1}$}  \vspace{-0.055em}
    \EndIf  \vspace{-0.055em}
    \vspace{-0.2em} 
    \State $\mathtt{v \leftarrow fs[m]}$ \Comment{larger value in next bump}  \vspace{-0.055em}
    \State $\mathtt{fs[m] \leftarrow m}$ \Comment{reset focus of $\mathtt{m}$}  \vspace{-0.055em}
\EndWhile  \vspace{-0.055em}
\vspace{-0.25em} 
\State \textbf{visit} $\mathtt{perm}$ \Comment{last $s$-word in Gray code}  \vspace{-0.055em}
\end{algorithmic}
\end{minipage}
\hfill
\begin{minipage}[t]{0.045\textwidth}
\raisebox{-4.8in}{\scalebox{-1}[1]{\includegraphics[page=6,angle=90,width=\textwidth,height=4.9in]{StirlingRopes.pdf}}}
\end{minipage}
\hfill
\begin{minipage}[t]{0.425\textwidth}
\vspace{-0.75em}
\small
\centering
\begin{tabular}[t]{|@{\;}c@{\;\;}c@{\;\;}c@{\;\;}c@{\;\;}c@{\;\;}c@{\;\;}c@{\;\;}c@{\;\;}c@{\;}|} \hline
\texttt{perm} & \texttt{v} & \texttt{u} & \texttt{i} & \texttt{j} & \texttt{left} & \texttt{inv} & \texttt{fs} & \texttt{dirs} \\ \hline
$\one\one\two\thr\thr\thr$ & $\thr$ & $\two$ & $\color{darkred}6$   & $\color{darkgreen}3$ & $\color{darkblue}1\color{darkgreen}3\color{darkred}4$ & $\color{darkblue}0\color{darkgreen}0\color{darkred}0$ & $123$ & $\color{darkblue}\texttt{-}\color{darkgreen}\texttt{-}\color{darkred}\texttt{-}$ \\[-0.1em]
$\one\one\thr\thr\thr\two$ & $\thr$ & $\one$ & $\color{darkred}5$   & $\color{darkblue}2$  & $\color{darkblue}1\color{darkgreen}6\color{darkred}3$ & $\color{darkblue}0\color{darkgreen}0\color{darkred}1$ & $123$ & $\color{darkblue}\texttt{-}\color{darkgreen}\texttt{-}\color{darkred}\texttt{-}$ \\[-0.1em]
$\one\thr\thr\thr\one\two$ & $\thr$ & $\one$ & $\color{darkred}4$   & $\color{darkblue}1$  & $\color{darkblue}1\color{darkgreen}6\color{darkred}2$ & $\color{darkblue}0\color{darkgreen}0\color{darkred}2$ & $123$ & $\color{darkblue}\texttt{-}\color{darkgreen}\texttt{-}\color{darkred}\texttt{-}$ \\[-0.1em]
$\thr\thr\thr\one\one\two$ & $\two$ & $\one$ & $\color{darkgreen}6$ & $\color{darkblue}5$  & $\color{darkblue}4\color{darkgreen}6\color{darkred}1$ & $\color{darkblue}0\color{darkgreen}0\color{darkred}3$ & $123$ & $\color{darkblue}\texttt{-}\color{darkgreen}\texttt{-}\color{darkred}\texttt{+}$ \\[-0.1em]
$\thr\thr\thr\one\two\one$ & $\thr$ & $\one$ & $\color{darkred}1$   & $\color{darkblue}4$  & $\color{darkblue}4\color{darkgreen}5\color{darkred}1$ & $\color{darkblue}0\color{darkgreen}1\color{darkred}3$ & $123$ & $\color{darkblue}\texttt{-}\color{darkgreen}\texttt{-}\color{darkred}\texttt{+}$ \\[-0.1em]
$\one\thr\thr\thr\two\one$ & $\thr$ & $\two$ & $\color{darkred}2$   & $\color{darkgreen}5$ & $\color{darkblue}1\color{darkgreen}5\color{darkred}2$ & $\color{darkblue}0\color{darkgreen}1\color{darkred}2$ & $123$ & $\color{darkblue}\texttt{-}\color{darkgreen}\texttt{-}\color{darkred}\texttt{+}$ \\[-0.1em]
$\one\two\thr\thr\thr\one$ & $\thr$ & $\one$ & $\color{darkred}3$   & $\color{darkblue}6$  & $\color{darkblue}1\color{darkgreen}2\color{darkred}3$ & $\color{darkblue}0\color{darkgreen}1\color{darkred}1$ & $123$ & $\color{darkblue}\texttt{-}\color{darkgreen}\texttt{-}\color{darkred}\texttt{+}$ \\[-0.1em]
$\one\two\one\thr\thr\thr$ & $\two$ & $\one$ & $\color{darkgreen}2$ & $\color{darkblue}1$  & $\color{darkblue}1\color{darkgreen}2\color{darkred}4$ & $\color{darkblue}0\color{darkgreen}1\color{darkred}0$ & $123$ & $\color{darkblue}\texttt{-}\color{darkgreen}\texttt{-}\color{darkred}\texttt{-}$ \\[-0.1em]
$\two\one\one\thr\thr\thr$ & $\thr$ & $\one$ & $\color{darkred}6$   & $\color{darkblue}3$  & $\color{darkblue}2\color{darkgreen}1\color{darkred}4$ & $\color{darkblue}0\color{darkgreen}2\color{darkred}0$ & $113$ & $\color{darkblue}\texttt{-}\color{darkgreen}\texttt{+}\color{darkred}\texttt{-}$ \\[-0.1em]
$\two\one\thr\thr\thr\one$ & $\thr$ & $\one$ & $\color{darkred}5$   & $\color{darkblue}2$  & $\color{darkblue}2\color{darkgreen}1\color{darkred}3$ & $\color{darkblue}0\color{darkgreen}2\color{darkred}1$ & $113$ & $\color{darkblue}\texttt{-}\color{darkgreen}\texttt{+}\color{darkred}\texttt{-}$ \\[-0.1em]
$\two\thr\thr\thr\one\one$ & $\thr$ & $\two$ & $\color{darkred}4$   & $\color{darkgreen}1$ & $\color{darkblue}5\color{darkgreen}1\color{darkred}2$ & $\color{darkblue}0\color{darkgreen}2\color{darkred}2$ & $113$ & $\color{darkblue}\texttt{-}\color{darkgreen}\texttt{+}\color{darkred}\texttt{-}$ \\[-0.1em]
$\thr\thr\thr\two\one\one$ & $\one$ &        &                      &                      & $\color{darkblue}5\color{darkgreen}4\color{darkred}1$ & $\color{darkblue}0\color{darkgreen}2\color{darkred}3$ & $123$ & $\color{darkblue}\texttt{-}\color{darkgreen}\texttt{+}\color{darkred}\texttt{+}$ \\[-0.1em]
\hline
\multicolumn{9}{c}{Variable trace at \textbf{visit} for $s=(2,1,3)$.} \\[-1.3em]
\multicolumn{9}{c}{}
\end{tabular}
\includegraphics[page=1,width=0.9\textwidth]{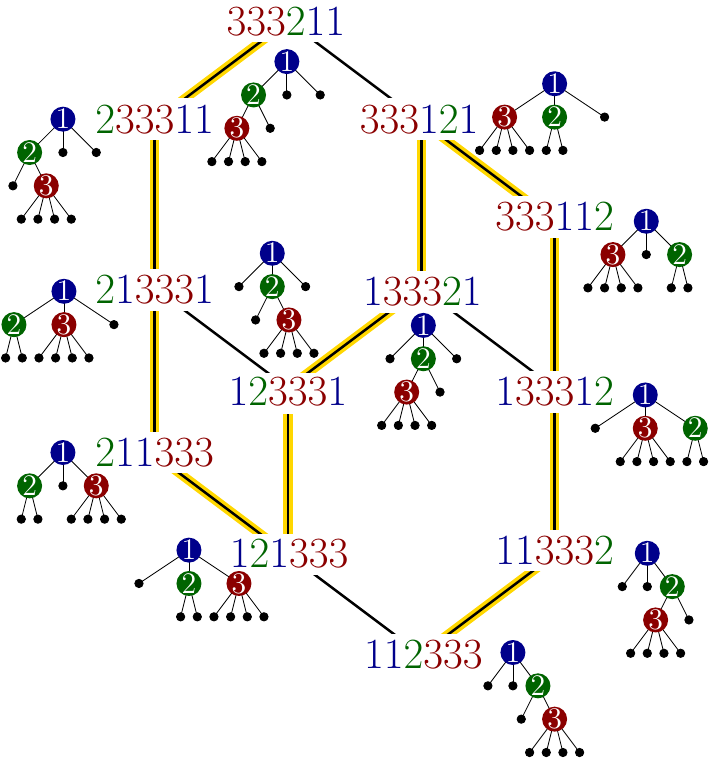}
\vspace{-1.1em}
\center{Figure 4: $s$-permutahedron Hamilton path for $s=(2,1,3)$ (or $s=(3,1,2)$ \cite{ceballos2019s}).} 
\end{minipage}
\end{algorithm}

\vspace{-2.6em}

\bibliographystyle{plain}
\bibliography{refs}

\end{document}